\documentclass[10pt,twoside]{article}

\usepackage[centertags]{amsmath}
\usepackage{amsfonts}
\usepackage{amssymb}
\usepackage{amsthm}
\usepackage{epsfig}
\usepackage{color}
\definecolor{c20}{rgb}{0.,0.7,0.}
\definecolor{c30}{rgb}{0.,0.,1.}
\definecolor{c40}{rgb}{1,0.1,0.7}
\definecolor{c50}{rgb}{1,0,0}
\definecolor{c60}{rgb}{1,0.9,0.1}

\allowdisplaybreaks[4]

\begin{document}
\baselineskip 15pt \setcounter{page}{1}
\title{\bf \Large  The asymptotic relation between the first crossing point and the last exit time of Gaussian order statistics sequences
\thanks{Research supported by Zhejiang Provincial Natural Science Foundation of China (Grant No. LY18A010020).} }
\author{{\small Zijun Ning;\ \ \ Zhongquan Tan\footnote{Corresponding author. E-mail address:  tzq728@163.com}}
\\
{\small\it  College of Data Science, Jiaxing University, Jiaxing 314001, PR China }\\
}
 \maketitle
 \baselineskip 15pt

\begin{quote}
{\bf Abstract:}\ \ In this paper, we study the asymptotic relation between the first crossing point and the last exit time for Gaussian order statistics which are generated by stationary weakly and strongly dependent Gaussian sequences.
It is shown that the first crossing point and the last exit time are asymptotically independent and dependent for weakly and strongly dependent respectively. The asymptotic relation between the first exit time and the last exit time for stationary weakly and strongly dependent Gaussian Gaussian sequences are also obtained.



\end{quote}

\section{Introduction}

The first crossing point (also known as the first exit time) and the last exit time  have been studied extensively in applied probability.
In risk theory, they play a very important role,  since the first exit time corresponds to the ruin time of a risk process and
the last exit time can be regarded as the final recovery time after when there will be no more ruin, see e.g., Gerber (1990).
In practice, when we consider the random noises problems, it usually adds the first zero crossing of sine wave to the random noises.
Instead of considering the zero crossing, we can transform it to stationary random sequence and
 consider therefore the first crossing of the sequence, see e.g., H\"{u}sler (1977).

Let $\{X_{i}, i\in \mathbb{Z}\}$ be a sequence of stationary random variables with $EX_{1}^{+}<\infty$, where $X^{+}$ means $\max\{X,0\}$. Define
$$\Upsilon(X)=\min\{i\in \mathbb{Z}: X_{i}>-\beta i\}$$ and
$$\tau(X)=\max\{i\geq0: X_{i}>\beta i\}$$
 the first crossing point and the last exit time of the sequence $\{X_{i}, i\in \mathbb{Z}\}$, respectively.
 Because of the existence of $EX_{1}^{+}$, both $\Upsilon(X)$ and $\tau(X)$ are well-defined and also finite almost surely, see e.g., H\"{u}sler (1979a).

H\"{u}sler (1977) first studied the limit distributions of first crossing point for stationary Gaussian sequences
and showed that with appropriately chosen normalization the limit distribution of the first crossing point and the maxima are the same.
H\"{u}sler (1979a) discussed the almost sure limiting behaviour including the a.s. stability and a.s. relative stability of the first crossing point for stationary Gaussian sequences and derived an iterated logarithm law for the first crossing point of the sequence.
It is easy to see from the definition of the first crossing point and the last exit time that the limit properties of them are equal.
The limit properties of last exit time for independently and identically distributed  random sequences and
stationary random sequences are derived by H\"{u}sler (1979b, 1980).  H\"{u}sler (1981) studied the law of the iterated logarithm for the last exit time of independent random sequences.
Recently, Karagodin and Lifshits (2021) considered a similar problem for continues time Gaussian processes and derived the limiting distribution of the scaled last exit time over a slowly growing linear boundary for stationary Gaussian processes.
All of the above mentioned results only consider the last exit time or the first crossing point for the linear boundary.
Karagodin (2022) studied the last exit time over a moving nonlinear boundary for a Gaussian process.
For some related results, see e.g., Burnashev and Golubev (2015) and Malinovskii (2018).

A natural question arisen is what is the joint limit distribution of the first and the last exit time.
For certain of Gaussian processes, H\"{u}sler and Zhang (2008) showed that
the conditional joint limit distribution of the first and the last exit time,
conditioned on ruin occurring, is a difference of two standard normal distributions.
Shao and Tan (2023) discussed the joint limit properties for the first crossing point and the last exit time for some depednent chi-sequences.

In this paper, we are interested in the asymptotic relation between the first crossing point and last exit time for Gaussian order statistics sequences. Let $\{X_{i}, i\in \mathbb{Z}\}$ be a sequence of stationary standard (mean 0, variance 1) Gaussian random variables  with  covariance function $r_{i}=E(X_{1}X_{i+1})$
and $\{X_{ij}, i\in \mathbb{Z}\}, j=1,2,\ldots,d$ with $d\geq 1$ be independent
copies of $\{X_{i}, i\in \mathbb{Z}\}$.
Define $\{O_{id}^{(r)}(X), i\in \mathbb{Z}\}$ a sequence of Gaussian $r$-th order statistics  generated by $X$ as follows
$$O_{id}^{(d)}(X):=\min_{j=1}^{d}X_{ij}\leq \cdots\leq O_{id}^{(r)}(X)\leq \cdots\leq O_{id}^{(1)}(X):=\max_{j=1}^{d}X_{ij},\ \ i\in \mathbb{Z},$$
where $r\in\{1,2,\ldots,d\}$.
The Gaussian $r$-th order statistics random variables play a very important role in applied fields,
for instance, in models concerned with the analysis of the surface roughness during
all machinery processes, see e.g., Alodat
(2011) and Worsley and Friston (2000). The limit properties of extremes of  Gaussian $r$-th order statistics sequences and processes are
studied in D\c{e}bicki et al. (2014, 2015, 2017), D\c{e}bicki and Kosinski (2018), Zhao (2018) and Song et al. (2022).

It is worth mentioning the following result, since it is highly relative to the main result of this paper.
Based on the results of D\c{e}bicki et al. (2017), Tan (2018) derived the weak limit theorems for the maxima of  Gaussian $r$-th order statistics variables.

\textbf{Theorem 1.1}. {\sl Let $\{O_{id}^{(r)}(X), i\in \mathbb{Z}\}$ be defined as above.
If the covariance function $r_{n}$  satisfies
$$r_{n}\ln n\rightarrow 0\ \ \mbox{as}\ \ n\rightarrow \infty,$$
then for any $x\in \mathbb{R}$
\begin{equation}
\label{Th1.1}
\lim_{n\rightarrow\infty}P\left(a_{n}\left(\max_{i=1}^{n}O_{id}^{(r)}(X)-b_{n}\right)\leq x\right)=\exp(-e^{-x}),
\end{equation}
where
\begin{equation}
\label{Th1.3}
a_{n}=(2r\ln n)^{1/2},\ \ b_{n}=\frac{1}{r}a_{n}+ a_{n}^{-1}\ln \Bigl(a_{n}^{-r} C_{d}^{r}(2\pi)^{-r/2}\Bigr)
\end{equation}
with $C_{d}^{r}=\frac{d!}{r!(d-r)!}$.
}

Theorem 1.1 obtained the limit distributions of maxima of Gaussian order statistics variables for weakly dependent cases. In the literature, assumption $r_{n}\ln n\rightarrow0 $ is referred to as
the weak dependence or the Berman's condition (see e.g. Leadbetter et al.(1983) and Piterbarg (1996)),
and, consequently, the stationary Gaussian sequences $\{X_{i}, i\in \mathbb{Z}\}$ is called a weakly dependent stationary
Gaussian sequence. We called also $\{O_{id}^{(r)}(X), i\in \mathbb{Z}\}$ weakly dependent Gaussian $r$-th order statistics sequence.
In analogy, $\{O_{id}^{(r)}(X), i\in \mathbb{Z}\}$  with correlation function satisfying assumption
$r_{n}\ln n\rightarrow\gamma>0 $ is called strongly dependent Gaussian $r$-th order statistics sequence.

In this paper, we will investigate the asymptotic relation between the first crossing point and the last exit time for Gaussian order statistics sequences. We showed that the first crossing point and the last exit time are asymptotic independent if the Gaussian order statistics sequences are weakly dependent, they are asymptotic dependent if the Gaussian order statistics sequences are strongly dependent.
We also showed that the limit properties of the first crossing point and the last exit time are very similar with that of the maxima.
Section 2 presents our main results and their proofs are provided in Section 3.

\section{Main results}

We define the first crossing point and the last exit time of Gaussian order statistics sequence $\{O_{id}^{(r)}(X), i\in \mathbb{Z}\}$
with function $i\beta$ as follows:
$$\Upsilon(O_{X}^{(r)})=\min\{i\in \mathbb{Z}: O_{id}^{(r)}(X)>-\beta i\}$$
and
$$\tau(O_{X}^{(r)})=\max\{i\geq 0: O_{id}^{(r)}(X)>\beta i\}$$
the last exit time of the sequence $\{X_{i}, i\in \mathbb{Z}\}$
with $\beta>0$. We consider the joint limit distribution of $\Upsilon_{O(X)}$ and $\tau(O_{X}^{(r)})$ as $\beta\rightarrow0$.

Now we state our main results.

\textbf{Theorem 2.1}. {\sl Let $\{O_{id}^{(r)}(X), i\in \mathbb{Z}\}$  be defined as above and define $C_{d}^{r}=\frac{d!}{r!(d-r!)}$.
If the covariance function $r_{n}$  satisfies
$$r_{n}\ln n\rightarrow 0\ \ \mbox{as}\ \ n\rightarrow \infty,$$
then for any $x,y\in \mathbb{R}$
\begin{equation}
\label{Th2.1.1}
\lim_{\beta\rightarrow0}P\left(\Upsilon(O_{X}^{(r)})\geq-u_{\beta}(x), \tau(O_{X}^{(r)})\leq u_{\beta}(y)\right)=\exp(-(e^{-x}+e^{-y}))
\end{equation}
and
\begin{equation}
\label{Th2.1.10}
\lim_{\beta\rightarrow0}P\left(\Upsilon(O_{X}^{(r)})\geq-u_{\beta}(x)\right)=\lim_{\beta\rightarrow0}P\left(\tau(O_{X}^{(r)})\leq u_{\beta}(x)\right)=\exp(-e^{-x}),
\end{equation}
where
\begin{equation}
\label{Th2.1.2}
u_{\beta}(x)=\left(a(\beta,r)+\frac{x}{ra(\beta,r)}\right)\beta^{-1}
\end{equation}
with
\begin{equation}
\label{Th2.1.3}
a(\beta,r)=\sqrt{\frac{2}{r}\left(\ln\beta^{-1}-\frac{r+1}{2}\ln\left(\frac{2}{r}\ln\beta^{-1}\right)+\ln(C_{d}^{r}r^{-1}(2\pi)^{-r/2})\right)}.
\end{equation}
}

For the strongly dependent case, we can only do the cases $r=1$ and $r=d$.

\textbf{Theorem 2.2}. {\sl Let $\{O_{id}^{(r)}(X), i\in \mathbb{Z}\}$  be defined as above.
Suppose that the covariance function $r_{n}$  satisfies
$$r_{n}\ln n\rightarrow \gamma>0\ \ \mbox{as}\ \ n\rightarrow \infty.$$
Define $\rho_{\beta}=\gamma/|\ln\beta|$ and $b(\beta,\rho_{\beta},r)=\sqrt{1-\rho_{\beta}}a(\beta/\sqrt{1-\rho_{\beta}},r)$ with $a(\cdot,\cdot)$ defined as in (\ref{Th2.1.3}) .\\
1.) For the case $r=1$, we have for any $x,y\in \mathbb{R}$
\begin{eqnarray}
\label{Th2.2.a}
&&\lim_{\beta\rightarrow0}P\left(\Upsilon(O_{X}^{(1)})\geq-u_{\beta}(x),\tau(O_{X}^{(1)})\leq u_{\beta}(y)\right)\nonumber\\
&&\ \ \ \ \ \ \ \ =\int_{\mathbb{R}^{d}}\exp\big(-\frac{1}{d}\sum_{j=1}^{d}(e^{-\sqrt{2\gamma}(x-z_{j})}+e^{-\sqrt{2\gamma}(y-z_{j})})\big)
d\Phi(z_{1})\cdots d\Phi(z_{d})
\end{eqnarray}
and
\begin{eqnarray}
\label{Th2.2.b}
\lim_{\beta\rightarrow0}P\left(\Upsilon(O_{X}^{(1)})\geq-u_{\beta}(x)\right)&=&\lim_{\beta\rightarrow0}P\left(\tau(O_{X}^{(1)})\leq u_{\beta}(x)\right)\nonumber\\
&=&\int_{\mathbb{R}^{d}}\exp\big(-\frac{1}{d}\sum_{j=1}^{d}e^{-\sqrt{2\gamma}(x-z_{j})}\big)
d\Phi(z_{1})\cdots d\Phi(z_{d}),
\end{eqnarray}
where
\begin{equation}
\label{Th2.2.c}
u_{\beta}(x)=\left(b(\beta,\rho_{\beta},1)+\sqrt{\rho_{\beta}}x\right)\beta^{-1}
\end{equation}
and $\Phi(\cdot)$ denotes the cumulative distribution function of a standard normal variable.\\
2.) For the case $r=d$, we have for any $x,y\in \mathbb{R}$
\begin{eqnarray}
\label{Th2.2.10}
&&\lim_{\beta\rightarrow0}P\left(\Upsilon(O_{X}^{(d)})\geq-u_{\beta}(x),\tau(O_{X}^{(d)})\leq u_{\beta}(y)\right)\nonumber\\
&&\ \ \ \ \ \ \ \ =\int_{\mathbb{R}^{d}}\exp\big(-(e^{-\sqrt{2d^{-1}\gamma}(x-\sum_{j=1}^{d}z_{j})}+e^{-\sqrt{2d^{-1}\gamma}(y-\sum_{j=1}^{d}z_{j})})\big)
d\Phi(z_{1})\cdots d\Phi(z_{d})
\end{eqnarray}
and
\begin{eqnarray}
\label{Th2.2.1}
\lim_{\beta\rightarrow0}P\left(\Upsilon(O_{X}^{(d)})\geq-u_{\beta}(x)\right)&=&\lim_{\beta\rightarrow0}P\left(\tau(O_{X}^{(d)})\leq u_{\beta}(x)\right)\nonumber\\
&=&\int_{\mathbb{R}^{d}}\exp\big(-e^{-\sqrt{2d^{-1}\gamma}(x-\sum_{j=1}^{d}z_{j})}\big)
d\Phi(z_{1})\cdots d\Phi(z_{d}),
\end{eqnarray}
where
\begin{equation}
\label{Th2.2.2}
u_{\beta}(x)=\left(b(\beta,\rho_{\beta},d)+\sqrt{\rho_{\beta}}d^{-1}x\right)\beta^{-1}.
\end{equation}
}

Next, let us consider the special case $d=1$, namely the stationary Gaussian case.

\textbf{Corollary 2.1}. {\sl Let $\{X_{i}, i\in \mathbb{Z}\}$ be a sequence of stationary standard
  Gaussian random variables  with  covariance function $r_{i}=E(X_{1}X_{i+1})$. \\
1.) Under the conditions of Theorem 2.1, we have
\begin{eqnarray*}
\label{Cor2.3.1}
\lim_{\beta\rightarrow0}P\left(\Upsilon(X)\geq -u_{\beta}(x), \tau(X)\leq u_{\beta}(y) \right)=\exp(-(e^{-x}+e^{-y})),
\end{eqnarray*}
where $u_{\beta}(x)$ is defined in (\ref{Th2.1.3}) with $r=d=1$.\\
2.) Under the conditions of Theorem 2.2, we have
\begin{eqnarray*}
\label{Cor2.3.2}
&&\lim_{\beta\rightarrow0}P\left(\Upsilon(X)\geq -u_{\beta}(x), \tau(X)\leq u_{\beta}(y) \right)\\
&&\ \ \ \ \ \ \ \ \ \ =\int_{\mathbb{R}}\exp(-(e^{-\sqrt{2\gamma}(x-z)}+e^{-\sqrt{2\gamma}(y-z)}))d\Phi(z),
\end{eqnarray*}
where $u_{\beta}(x)$ is defined in (\ref{Th2.2.2}) with $d=1$.
}

Corollary 2.1 extends Theorems 2.1 and 2.2 of H\"{u}sler (1977) which derived only the limit distribution for $\Upsilon(X)$.  But there is still a case which can not be derived
from Theorems 2.1 and 2.2 and will be discussed in next theorem.

\textbf{Theorem 2.3}. {\sl Let $\{X_{i}, i\in \mathbb{Z}\}$ be a sequence of stationary standard
  Gaussian random variables  with  covariance function $r_{i}=E(X_{1}X_{i+1})$.
Suppose that $r_{n}$ is convex for all $n\geq 0$, $r_{n}=o(1)$, and $(r_{n}\ln n)^{-1}$
is monotone for large $n$ and $o(1)$.
Then for any $x\in \mathbb{R}$
\begin{equation}
\label{Th2.3.1}
\lim_{\beta\rightarrow0}P\left(\Upsilon(X)\geq -u_{\beta}(x), \tau(X)\leq u_{\beta}(y) \right)=\Phi(\min{x,y}),
\end{equation}
where $u_{\beta}(x)$ is defined in (\ref{Th2.2.2}) with $d=1$ and $\rho_{\beta}=r_{[\beta^{-1}]}$. Here $[x]$ means the integral parts of $x$.
}

\textbf{Remark 2.1}. i). It can be seen that with appropriately chosen normalization the limit distributions of $\max_{i=1}^{n}O_{id}^{(r)}(X)$ $(n\rightarrow\infty)$, $\tau(O_{X}^{(r)})$ and $\Upsilon(O_{X}^{(r)})$ $(\beta\rightarrow0)$ are the same for the weakly dependent case.\\
ii). As in H\"{u}sler (1977), the above results can be generalized for the first crossing point $\Upsilon(O_{X}^{(r)})=\min\{i\in \mathbb{Z}: O_{id}^{(r)}(X)>-\beta f(i)\}$ and the last exit time $\tau(O_{X}^{(r)})=\max\{i\geq 0: O_{id}^{(r)}(X)>\beta f(i)\}$ with $f(t)=t^{\vartheta}L(t)$ for $t\rightarrow\infty$, where $\vartheta>0$ and $L(t)$ varies slowly.\\

\section{Proofs }

In this section, we give the proofs.
As usual, $a_{n}\ll b_{n}$ means $\limsup_{n\rightarrow\infty}|a_{n}/b_{n}|<+\infty$.
$K$ will denote a constant whose value will change from line to line.

The following lemma plays a crucial role in our proofs.

\textbf{Lemma 3.1}. {\sl
Denote by
$\mathcal{Y}=(Y_{il})_{n\times d}$ and $\mathcal{Z}=(Z_{il})_{n\times d}$  two random arrays with $N(0,1)$ components, and let
$(\sigma^{(1)}_{il,jk})_{{dn \times dn}}$ and $(\sigma^{(0)}_{il,jk})_{dn \times dn}$ be the covariance matrices of $\mathcal{Y}$ and $\mathcal{Z}$, respectively, with
$ \sigma^{(1)}_{il,jk} := E{Y_{il}Y_{jk}}$ and $\sigma^{(0)}_{il,jk} := E{Z_{il}Z_{jk}}, 1\le i,j\le n, 1\le l,k\le d.$
Suppose that the columns of both $\mathcal{Y}$ and $\mathcal{Z}$ are mutually independent, i.e.,
$$\sigma^{(\kappa)}_{il,jk}=\sigma^{(\kappa)}_{ik,jk}\mathbf{1}(k=l), 1\leq i,j\leq n, 1\leq k,l\leq d, \kappa=0,1.$$
Further, suppose that $\sigma_{i,j}=\max_{1\leq k,l\leq d}\{|\sigma^{(0)}_{ik,jk}|,|\sigma^{(1)}_{ik,jk}|\}<1$,
where $\mathbf{1}(\cdot)$ will denote an  indicator function.
Define $(O_{1d}^{(r)}(Y), \cdots, O_{nd}^{(r)}(Y))$ to be the  $r$-th order statistics vector generated by $\mathcal{X}$ as follows
\begin{eqnarray*}
\min_{1\le k\le d}Y_{ik}= O_{id}^{(d)}(Y) \leq  \cdots \leq O_{id}^{(r)}(Y)\leq\cdots \leq O_{id}^{(1)}(Y)= \max_{1\le k\le d} Y_{ik},\quad 1\le i\le n.
\end{eqnarray*}
Similarly, we write  $(O_{1d}^{(r)}(Z), \cdots, O_{nd}^{(r)}(Z))$ which is generated by $\mathcal{Z}$. Then for any real numbers $u_{1},\ldots,u_{n}$,
and any $1\leq r\leq d$
\begin{eqnarray*}
&&\bigg|P\bigg(O_{id}^{(r)}(Y)\leq u_{i},i=1,\ldots,n\bigg)-P\bigg(O_{id}^{(r)}(Z)\leq u_{i},i=1,\ldots,n\bigg)\bigg|\\
&&\leq K\sum_{k=1}^{d}\sum_{1\leq i< j\leq n}\frac{|\sigma^{(1)}_{ik,jk}-\sigma^{(0)}_{ik,jk}|}{(u_{i}u_{j})^{r-1}}\exp\left(-\frac{r(u_{i}^{2}+u_{j}^{2})}{2(1+\sigma_{i,j})}\right).
\end{eqnarray*}
}

\textbf{Proof}.
Let $(Z_i, Z_j)$ be a bivariate standard normal random vector with correlation $|\delta_{ij}|$.  By a
similar argument as the proof on p. 225 of Leadbetter et al. (1983), we can show for $u_{j}\geq u_{i}$
\begin{eqnarray*}
&&P\{Z_i >u_{i}, Z_j >u_{j}\} \leq \frac{2(1-\delta_{ij}^{2})}{u_{i}(u_{j}-|\delta_{ij}|u_{i})}\phi(u_{i}, u_{j}; |\delta_{ij}|).
\end{eqnarray*}
where $\phi(u,v; \delta)$ is the probability density function of a two dimensional normal random variables. If $u_{i}=O( u_{j})$,
we have
$$\frac{2(1-\delta_{ij}^{2})}{u_{i}(u_{j}-|\delta_{ij}|u_{i})}\leq \frac{2(1-\delta_{ij}^{2})}{u_{i}u_{j}|1-O(1)|\delta_{ij}||}\leq \frac{K}{u_{i}u_{j}}.$$
If $u_{i}=o( u_{j})$, obviously,
$$\frac{2(1-\delta_{ij}^{2})}{u_{i}(u_{j}-|\delta_{ij}|u_{i})}= O(1)\frac{2(1-\delta_{ij}^{2})}{u_{i}u_{j}}\leq \frac{K}{u_{i}u_{j}}.$$
By the same arguments, we can get the same bound for the case $u_{j}\leq u_{i}$. Thus, we have
\begin{eqnarray*}
&&P\{Z_i >u_{i}, Z_j >u_{j}\} \leq \frac{K}{u_{i}u_{j}}\phi(u_{i}, u_{j}; |\delta_{ij}|).
\end{eqnarray*}
Using the above
inequality to replace (4.28) in the proof of Theorem 2.4 of D\c{e}bicki et al. (2017), we can prove the lemma.

\subsection{Proof of Theorem 2.1}

We need the following lemma to prove Theorem 2.1.

\textbf{Lemma 3.2}. {\sl Let $\{\xi_{i}, i\in \mathbb{Z}\}$ be a sequence of independent standard Gaussian random variable and $\{O_{id}^{(r)}(\xi), i\in \mathbb{Z}\}$ be a sequence of Gaussian $r$-th order statistics generated by $\xi$.
We have for any $x,y\in \mathbb{R}$
\begin{equation}
\label{Le3.2.1}
\lim_{\beta\rightarrow0}P\left(\Upsilon(O_{\xi}^{(r)})\geq-u_{\beta}(x),\tau(O_{\xi}^{(r)})\leq u_{\beta}(y)\right)=\exp(-(e^{-x}+e^{-y})),
\end{equation}
wehre $u_{\beta}(x)$ is defined in  (\ref{Th2.1.2}).
}

\textbf{Proof}. First note that $\beta u_{\beta}(x)\rightarrow\infty$ as $\beta\rightarrow0$.
From Lemma 2 of D\c{e}bicki et al. (2015), we have
\begin{equation}
\label{Le3.2.2}
P(O_{1d}^{(r)}(\xi)\geq i\beta)=C_{d}^{r}(1-\Phi(i\beta))^{r}(1+o(1))
\end{equation}
as $i\beta\rightarrow\infty$.
Since $O_{id}^{(r)}(\xi)$ are independent, using the following well-known fact
\begin{equation}
\label{Le3.2.3}
1-\Phi(x)=\frac{1}{\sqrt{2\pi}x}e^{-\frac{x^{2}}{2}}(1+o(1))
\end{equation}
as $x\rightarrow\infty$, we have as $\beta\rightarrow0$
\begin{eqnarray*}
-\ln P\left(\Upsilon(O_{\xi}^{(r)})\geq -u_{\beta}(x)\right)
&=&-\ln \prod_{i=u_{\beta}(x)}^{\infty}P(O_{1d}^{(r)}(\xi)\leq i\beta)\\
&=&-\sum_{i=u_{\beta}(x)}^{\infty}\ln (1- C_{d}^{r}(1-\Phi(i\beta)^{r})(1+o(1)))\\
&=&C_{d}^{r}(2\pi)^{-r/2}r^{-1}\beta^{-1}(\beta u_{\beta}(x))^{-(r+1)}e^{-\frac{1}{2}r(\beta u_{\beta}(x))^{2}}(1+o(1)).
\end{eqnarray*}
Since $(\beta u_{\beta}(x))^{2}=(a(\beta,r))^{2}+2x/r+o(1)$, in view of (\ref{Th2.1.3}), we have
$$-\ln P\left(\Upsilon(O_{\xi}^{(r)})\geq -u_{\beta}(x)\right)=e^{-x}(1+o(1))$$
as $\beta\rightarrow0$.
Similarly, we have
$$-\ln P\left(\tau(O_{\xi}^{(r)})\leq u_{\beta}(y)\right)=e^{-y}(1+o(1))$$
as $\beta\rightarrow0$.
By the independence of $\{O_{id}^{(r)}(\xi), i\in \mathbb{Z}\}$, we can see that
\begin{eqnarray*}
&&P\left(\Upsilon(O_{\xi}^{(r)})\geq-u_{\beta}(x),\tau(O_{\xi}^{(r)})\leq u_{\beta}(y)\right)\\
&&=P\left(O_{jd}^{(r)}(\xi)\leq -j\beta, j\leq -u_{\beta}(x),O_{id}^{(r)}(\xi)\leq i\beta, i\geq u_{\beta}(y)\right)\\
&&=P\left(O_{jd}^{(r)}(\xi)\leq -j\beta, j\leq -u_{\beta}(x)\right)P\left(O_{id}^{(r)}(\xi)\leq i\beta, i\geq u_{\beta}(y)\right)
\end{eqnarray*}
and therefore the lemma follows.

\textbf{Proof of Theorem 2.1}. By Lemma 3.2, to prove Theorem 2.1, we only need to show that
$$|P\left(\Upsilon(O_{X}^{(r)})\geq-u_{\beta}(x), \tau(O_{X}^{(r)})\leq u_{\beta}(y)\right)-P\left(\Upsilon(O_{\xi}^{(r)})\geq-u_{\beta}(x), \tau(O_{\xi}^{(r)})\leq u_{\beta}(y)\right)|\rightarrow0$$
as $\beta\rightarrow0$. We use Lemma 3.1 to bound this difference by
\begin{eqnarray*}
&&K\sum_{i=1}^{\infty}|r_{i}|\sum_{j=u_{\beta}(x)}^{\infty}(j\beta)^{-2(r-1)}\exp\left(-\frac{r((j\beta)^{2}+((j+i)\beta)^{2})}{2(1+|r_{i}|)}\right)\\
&&+K\sum_{i=1}^{\infty}|r_{i}|\sum_{j=u_{\beta}(y)}^{\infty}(j\beta)^{-2(r-1)}\exp\left(-\frac{r((j\beta)^{2}+((j+i)\beta)^{2})}{2(1+|r_{i}|)}\right)\\
&&+K\sum_{i=u_{\beta}(x)}^{\infty}\sum_{j=-\infty}^{-u_{\beta}(y)}|r_{|i-j|}|(i\beta)^{-(r-1)}(-j\beta)^{-(r-1)}
\exp\left(-\frac{r((i\beta)^{2}+(j\beta)^{2})}{2(1+|r_{|i-j|}|)}\right)\\
&&=:K(I_{\beta}+J_{\beta}+K_{\beta}).
\end{eqnarray*}
Define $\delta=\sup_{i\geq 1}|r_{i}|<1$ and $\delta_{n}=\sup_{i\geq n}|r_{i}|<1$. Choose $\varepsilon$ such that
$\varepsilon<(1-\delta)/(1+\delta)$. Split the above sum $I_{\beta}$ into two parts as
\begin{eqnarray*}
I_{\beta}&=&\sum_{i=1}^{[\beta^{-\varepsilon}]}|r_{i}|\sum_{j=u_{\beta}(x)}^{\infty}(j\beta)^{-2(r-1)}
\exp\left(-\frac{r((j\beta)^{2}+((j+i)\beta)^{2})}{2(1+|r_{i}|)}\right)\\
&+&\sum_{i=[\beta^{-\varepsilon}]+1}^{\infty}|r_{i}|\sum_{j=u_{\beta}(x)}^{\infty}(j\beta)^{-2(r-1)}
\exp\left(-\frac{r((j\beta)^{2}+((j+i)\beta)^{2})}{2(1+|r_{i}|)}\right)\\
&=:&I_{\beta,1}+I_{\beta,2}.
\end{eqnarray*}
For the first term $I_{\beta,1}$, using (\ref{Th2.1.2}) and (\ref{Th2.1.3}), we have
\begin{eqnarray*}
I_{\beta,1}&\leq &\beta^{-\varepsilon}\sum_{j=u_{\beta}(x)}^{\infty}(j\beta)^{-2(r-1)}
\exp\left(-\frac{r(j\beta)^{2}}{1+\delta}\right)\\
&\ll& \beta^{-(\varepsilon+1)}(u_{\beta}(x)\beta)^{-(2r-1)}\exp\left(-\frac{r(u_{\beta}(x)\beta)^{2}}{1+\delta}\right)\\
&\ll & \beta^{-\varepsilon+\frac{1-\delta}{1+\delta}}|\ln\beta|^{\frac{1}{2}+\frac{1-r\delta}{1+\delta}},
\end{eqnarray*}
which tends to 0 as $\beta\rightarrow0$ by the choice of $\varepsilon$.
For the term $I_{\beta,2}$, recall that $\delta_{[\beta^{-\varepsilon}]}=\sup_{i\geq [\beta^{-\varepsilon}]}|r_{i}|$ and note that $r_{n}\ln n\rightarrow 0$ as $n\rightarrow\infty$ implies that $\delta_{[\beta^{-\varepsilon}]}|\ln \beta|\rightarrow 0$
as $\beta\rightarrow0$.  We have in view of (\ref{Th2.1.2}) and (\ref{Th2.1.3}) again
\begin{eqnarray*}
I_{\beta,2}&\leq &\delta_{[\beta^{-\varepsilon}]} \sum_{i=[\beta^{-\varepsilon}]}^{\infty}\exp\left(-\frac{r((i+u_{\beta}(x))\beta)^{2}}{2(1+\delta_{[\beta^{-\varepsilon}]})}\right)
\sum_{j=u_{\beta}(x)}^{\infty}(j\beta)^{-2(r-1)}\exp\left(-\frac{r(j\beta)^{2}}{2(1+\delta_{[\beta^{-\varepsilon}]})}\right)\\
&\ll&\delta_{[\beta^{-\varepsilon}]} \beta^{-2}(\beta u_{\beta}(x))^{-2r} \exp\left(-\frac{r(u_{\beta}(x)\beta)^{2}}{1+\delta_{[\beta^{-\varepsilon}]}}\right)\\
&\ll & \delta_{[\beta^{-\varepsilon}]}\beta^{-2\delta_{[\beta^{-\varepsilon}]}}|\ln\beta|^{1-(1+r)\delta_{[\beta^{-\varepsilon}]}},
\end{eqnarray*}
which tends to 0 as $\beta\rightarrow0$ since $\delta_{[\beta^{-\varepsilon}]}|\ln \beta|\rightarrow0$ as $\beta\rightarrow0$.
Thus, $I_{\beta}\rightarrow\infty$ as $\beta\rightarrow0$. Similarly, we can show that $J_{\beta}\rightarrow\infty$ as $\beta\rightarrow0$.

For the term $K_{\beta}$,
note that $u_{\beta}(x)\sim u_{\beta}(y)\sim u_{\beta}:=\beta^{-1}(\frac{2}{r}\ln \beta^{-1})^{1/2}$ as $\beta\rightarrow0$.
Thus, by the same arguments as for $I_{\beta,2}$, we have
\begin{eqnarray*}
K_{\beta}&\leq &\delta_{[u_{\beta}]} \sum_{i=[u_{\beta}(x)]}^{\infty}(i\beta)^{-(r-1)}\exp\left(-\frac{r(i\beta)^{2}}{2(1+\delta_{[u_{\beta}]})}\right)
\sum_{j=[u_{\beta}(y)]}^{\infty}(j\beta)^{-(r-1)}\exp\left(-\frac{r(j\beta)^{2}}{2(1+\delta_{[u_{\beta}]})}\right)\\
&\ll&\delta_{[u_{\beta}]} \beta^{-2}(\beta u_{\beta})^{-2r} \exp\left(-\frac{r(u_{\beta}(x)\beta)^{2}}{2(1+\delta_{[u_{\beta}]})}\right)\exp\left(-\frac{r(u_{\beta}(y)\beta)^{2}}{2(1+\delta_{[u_{\beta}]})}\right)\\
&\ll & \delta_{[u_{\beta}]}\beta^{-2\delta_{[u_{\beta}]}}|\ln\beta|^{1-(1+r)\delta_{[u_{\beta}]}},
\end{eqnarray*}
which tends to 0 as $\beta\rightarrow0$ since $\delta_{[u_{\beta}]}|\ln \beta|\rightarrow0$ as $\beta\rightarrow0$.
 Thus, we have
$K_{\beta}\rightarrow0$  as $\beta\rightarrow0$.
The proof of Theorem 2.1 is complete.

\subsection{Proof of Theorem 2.2}

To prove Theorem 2.2, we need the following two lemmas.

\textbf{Lemma 3.3}. {\sl Let $\rho_{\beta}$ be  a function of $\beta$ with $0< \rho_{\beta}<c<1$ and
$Y_{ij}=\sqrt{1-\rho_{\beta}}\xi_{ij}+\sqrt{\rho_{\beta}}U_{j}$, $i\in \mathbb{N}, j=1,2\ldots,d$, where $\xi_{ij}$ are independent Gaussian random
vectors and $U_{j}, j=1,2\ldots,d$ are a standard normal random variables independent with $\xi_{ij}$. Denote by $O_{id}^{(r)}(Y)$ the Gaussian order statistics sequences generated by $Y_{ij}$.
Suppose that $\rho_{\beta}|\ln \beta|\rightarrow\gamma>0$ as $\beta\rightarrow0$.\\
1). If $r=1$, we have for any $x,y\in \mathbb{R}$
\begin{eqnarray*}
\label{Le3.4.0}
&&\lim_{\beta\rightarrow0}P\left(O_{id}^{(1)}(Y)\leq -i\beta, -u_{\beta}(x)-\beta^{-1}<i\leq -u_{\beta}(x), O_{kd}^{(1)}(Y)\leq k\beta, u_{\beta}(y)\leq k< u_{\beta}(y)+\beta^{-1} \right)\nonumber\\
&&\ \ \ \ \ \ \ \ \ =\int_{\mathbb{R}^{d}}\exp\big(-\frac{1}{d}\sum_{j=1}^{d}(e^{-\sqrt{2\gamma}(x-z_{j})}+e^{-\sqrt{2\gamma}(y-z_{j})})\big)
d\Phi(z_{1})\cdots d\Phi(z_{d}),
\end{eqnarray*}
where $u_{\beta}(x)$ is defined in (\ref{Th2.2.c}) with $\rho_{\beta}$ satisfying the above assumptions.\\
2). If $r=d$, we have for any $x,y\in \mathbb{R}$
\begin{eqnarray*}
\label{Le3.4.1}
&&\lim_{\beta\rightarrow0}P\left(O_{id}^{(d)}(Y)\leq -i\beta, -u_{\beta}(x)-\beta^{-1}<i\leq -u_{\beta}(x), O_{kd}^{(d)}(Y)\leq k\beta, u_{\beta}(y)\leq k< u_{\beta}(y)+\beta^{-1} \right)\nonumber\\
&&\ \ \ \ \ \ \ \ \ =\int_{\mathbb{R}^{d}}\exp\{-(e^{\sqrt{2d^{-1}\gamma}(x-\sum_{j=1}^{d}z_{j})}+e^{\sqrt{2d^{-1}\gamma}(y-\sum_{j=1}^{d}z_{j})})\}d\Phi(z_{1})\cdots d\Phi(z_{d}),
\end{eqnarray*}
where $u_{\beta}(x)$ is defined in (\ref{Th2.2.2}) with $\rho_{\beta}$ satisfying the above assumptions.
}

\textbf{Proof}. We show first the case $r=1$. By the definition of $Y_{ij}$ and the independence of $\xi_{ij}$, we have
\begin{eqnarray}
\label{Le3.4.3}
&&P\left(O_{id}^{(1)}(Y)\leq -i\beta, -u_{\beta}(x)-\beta^{-1}<i\leq -u_{\beta}(x),O_{kd}^{(1)}(Y)\leq k\beta, u_{\beta}(y)\leq k< u_{\beta}(y)+\beta^{-1}  \right)\nonumber\\
&&=P\left(\max_{j=1}^{d}\{\sqrt{1-\rho_{\beta}}\xi_{ij}+\sqrt{\rho_{\beta}}U_{j}\}\leq -i\beta, -u_{\beta}(x)-\beta^{-1}<i\leq -u_{\beta}(x)\right.\nonumber\\
&&\ \ \ \ \ \ \ \ \ \ \ \left. \max_{j=1}^{d}\{\sqrt{1-\rho_{\beta}}\xi_{kj}+\sqrt{\rho_{\beta}}U_{j}\}\leq k\beta, u_{\beta}(y)\leq k< u_{\beta}(y)+\beta^{-1} \right)\nonumber\\
&&=\int_{\mathbb{R}^{d}}
P\left(\max_{j=1}^{d}\{\sqrt{1-\rho_{\beta}}\xi_{ij}+\sqrt{\rho_{\beta}}z_{j}\}\leq -i\beta,  -u_{\beta}(x)-\beta^{-1}<i\leq -u_{\beta}(x),\right. \nonumber\\
&&\ \ \ \ \ \ \ \ \ \ \ \left. \max_{j=1}^{d}\{\sqrt{1-\rho_{\beta}}\xi_{kj}+\sqrt{\rho_{\beta}}z_{j}\}\leq k\beta, u_{\beta}(y)\leq k< u_{\beta}(y)+\beta^{-1}  \right)d\Phi(z_{1})\cdots d\Phi(z_{d})\nonumber\\
&&=\int_{\mathbb{R}^{d}}
\prod_{i=-u_{\beta}(x)-[\beta^{-1}]+1}^{-u_{\beta}(x)}\prod_{j=1}^{d}P\left(\sqrt{1-\rho_{\beta}}\xi_{ij}+\sqrt{\rho_{\beta}}z_{j}\leq -i\beta \right)\nonumber\\
&&\ \ \ \ \ \ \ \ \ \ \times\prod_{k=u_{\beta}(y)}^{u_{\beta}(y)+[\beta^{-1}]+1}\prod_{j=1}^{d}P\left(\sqrt{1-\rho_{\beta}}\xi_{kj}+\sqrt{\rho_{\beta}}z_{j}\leq k\beta \right) d\Phi(z_{1})\cdots d\Phi(z_{d})\nonumber\\
&&=:\int_{\mathbb{R}^{d}}H(\beta,x,z_{1},\ldots,z_{d})H(\beta,y,z_{1},\ldots,z_{d})d\Phi(z_{1})\cdots d\Phi(z_{d}).
\end{eqnarray}
By the same arguments as for the proof of Lemma 3.2, we have
\begin{eqnarray*}
&&-\ln H(\beta,x,z_{1},\ldots,z_{d})\\
&&=-\ln \prod_{i=-u_{\beta}(x)-[\beta^{-1}]+1}^{-u_{\beta}(x)}\prod_{j=1}^{d}\left[1-P\left(\xi_{ij}\geq \frac{-(i\beta+\sqrt{\rho_{\beta}}z_{j})}{\sqrt{1-\rho_{\beta}}} \right)\right]\\
&&=-\sum_{i=-u_{\beta}(x)-[\beta^{-1}]+1}^{-u_{\beta}(x)}\sum_{j=1}^{d}\ln \left[1-P\left(\xi_{ij}\geq \frac{-(i\beta+\sqrt{\rho_{\beta}}z_{j})}{\sqrt{1-\rho_{\beta}}} \right)\right]\\
&&=\sum_{i=-u_{\beta}(x)-[\beta^{-1}]+1}^{-u_{\beta}(x)}\sum_{j=1}^{d}
\left(1-\Phi\left(\frac{-(i\beta+\sqrt{\rho_{\beta}}z_{j})}{\sqrt{1-\rho_{\beta}}}\right)\right)(1+o(1))\\
&&=\sum_{i=-u_{\beta}(x)-[\beta^{-1}]+1}^{-u_{\beta}(x)}\sum_{j=1}^{d}(2\pi)^{-1/2}(1-\rho_{\beta})^{1/2}|i\beta|^{-1}
\exp\left(-\frac{(i\beta)^{2}+2i\beta\sqrt{\rho_{\beta}}z_{j}}{2(1-\rho_{\beta})}\right)(1+o(1))\\
&&=(2\pi)^{-1/2}(1-\rho_{\beta})^{1/2}\beta^{-1}(\beta u_{\beta}(x))^{-2}
\exp\left(-\frac{(\beta u_{\beta}(x))^{2}}{2(1-\rho_{\beta})}\right)
\sum_{j=1}^{d}\exp\left(\frac{\beta u_{\beta}(x)\sqrt{\rho_{\beta}}z_{j}}{1-\rho_{\beta}}\right)(1+o(1)).
\end{eqnarray*}
Now, using (\ref{Th2.2.c}) and the fact that $\rho_{\beta}|\ln \beta|\rightarrow\gamma>0$ as $\beta\rightarrow0$, we have
\begin{eqnarray}
\label{Le3.4.4}
-\ln H(\beta,x,z_{1},\ldots,z_{d})\rightarrow d^{-1}\sum_{j=1}^{d}\exp\big(-\sqrt{2\gamma}(x-z_{j})\big),
\end{eqnarray}
as $\beta\rightarrow0$. Similarly, we have
\begin{eqnarray}
\label{Le3.4.5}
-\ln H(\beta,y,z_{1},\ldots,z_{d})\rightarrow  d^{-1}\sum_{j=1}^{d}\exp\big(-\sqrt{2\gamma}(y-z_{j})\big),
\end{eqnarray}
as $\beta\rightarrow0$.
Now, (\ref{Le3.4.3}-\ref{Le3.4.5}) combining with the dominated convergence theorem completes the proof of the case $r=1$.\\
Next, we prove the case $r=d$. By the definition of $Y_{ij}$ and the independence of $\xi_{ij}$ again, we have
\begin{eqnarray}
\label{Le3.4.30}
&&P\left(O_{id}^{(d)}(Y)\leq -i\beta, -u_{\beta}(x)-\beta^{-1}<i\leq -u_{\beta}(x),O_{kd}^{(d)}(Y)\leq k\beta, u_{\beta}(y)\leq k< u_{\beta}(y)+\beta^{-1}  \right)\nonumber\\
&&=P\left(\min_{j=1}^{d}\{\sqrt{1-\rho_{\beta}}\xi_{ij}+\sqrt{\rho_{\beta}}U_{j}\}\leq -i\beta, -u_{\beta}(x)-\beta^{-1}<i\leq -u_{\beta}(x)\right.\nonumber\\
&&\ \ \ \ \ \ \ \ \ \ \ \left. \min_{j=1}^{d}\{\sqrt{1-\rho_{\beta}}\xi_{kj}+\sqrt{\rho_{\beta}}U_{j}\}\leq k\beta, u_{\beta}(y)\leq k< u_{\beta}(y)+\beta^{-1} \right)\nonumber\\
&&=\int_{\mathbb{R}^{d}}
P\left(\min_{j=1}^{d}\{\sqrt{1-\rho_{\beta}}\xi_{ij}+\sqrt{\rho_{\beta}}z_{j}\}\leq -i\beta,  -u_{\beta}(x)-\beta^{-1}<i\leq -u_{\beta}(x),\right. \nonumber\\
&&\ \ \ \ \ \ \ \ \ \ \ \left. \min_{j=1}^{d}\{\sqrt{1-\rho_{\beta}}\xi_{kj}+\sqrt{\rho_{\beta}}z_{j}\}\leq k\beta, u_{\beta}(y)\leq k< u_{\beta}(y)+\beta^{-1}  \right)d\Phi(z_{1})\cdots d\Phi(z_{d})\nonumber\\
&&=\int_{\mathbb{R}^{d}}
\prod_{i=-u_{\beta}(x)-[\beta^{-1}]+1}^{-u_{\beta}(x)}P\left(\min_{j=1}^{d}\{\sqrt{1-\rho_{\beta}}\xi_{ij}+\sqrt{\rho_{\beta}}z_{j}\}\leq -i\beta \right)\nonumber\\
&&\ \ \ \ \ \ \ \ \ \ \times\prod_{k=u_{\beta}(y)}^{u_{\beta}(y)+[\beta^{-1}]+1}P\left(\min_{j=1}^{d}\{\sqrt{1-\rho_{\beta}}\xi_{kj}+\sqrt{\rho_{\beta}}z_{j}\}\leq k\beta \right) d\Phi(z_{1})\cdots d\Phi(z_{d})\nonumber\\
&&=:\int_{\mathbb{R}^{d}}G(\beta,x,z_{1},\ldots,z_{d})G(\beta,y,z_{1},\ldots,z_{d})d\Phi(z_{1})\cdots d\Phi(z_{d}).
\end{eqnarray}
By the same arguments as for the proof of Lemma 3.2, we have
\begin{eqnarray*}
&&-\ln G(\beta,x,z_{1},\ldots,z_{d})\\
&&=-\ln \prod_{i=-u_{\beta}(x)-[\beta^{-1}]+1}^{-u_{\beta}(x)}\left[1-\prod_{j=1}^{d}P\left(\xi_{ij}\geq \frac{-(i\beta+\sqrt{\rho_{\beta}}z_{j})}{\sqrt{1-\rho_{\beta}}} \right)\right]\\
&&=-\sum_{i=-u_{\beta}(x)-[\beta^{-1}]+1}^{-u_{\beta}(x)}\ln \left[1-\prod_{j=1}^{d}P\left(\xi_{ij}\geq \frac{-(i\beta+\sqrt{\rho_{\beta}}z_{j})}{\sqrt{1-\rho_{\beta}}} \right)\right]\\
&&=\sum_{i=-u_{\beta}(x)-[\beta^{-1}]+1}^{-u_{\beta}(x)}\prod_{j=1}^{d}
\left(1-\Phi\left(\frac{-(i\beta+\sqrt{\rho_{\beta}}z_{j})}{\sqrt{1-\rho_{\beta}}}\right)\right)(1+o(1))\\
&&=\sum_{i=-u_{\beta}(x)-[\beta^{-1}]+1}^{-u_{\beta}(x)}(2\pi)^{-d/2}(1-\rho_{\beta})^{d/2}|i\beta|^{-d}
\exp\left(-\frac{d(i\beta)^{2}+2i\beta\sqrt{\rho_{\beta}}\sum_{j=1}^{d}z_{j}}{2(1-\rho_{\beta})}\right)(1+o(1))\\
&&=(2\pi)^{-d/2}d^{-1}(1-\rho_{\beta})^{d/2}\beta^{-1}(\beta u_{\beta}(x))^{-(d+1)}
\exp\left(-\frac{d(\beta u_{\beta}(x))^{2}}{2(1-\rho_{\beta})}\right)
\exp\left(\frac{\beta u_{\beta}(x)\sqrt{\rho_{\beta}}\sum_{j=1}^{d}z_{j}}{1-\rho_{\beta}}\right)(1+o(1)).
\end{eqnarray*}
Now, using (\ref{Th2.2.2}) and the fact $\rho_{\beta}|\ln \beta|\rightarrow\gamma>0$ as $\beta\rightarrow0$, we have
\begin{eqnarray}
\label{Le3.4.40}
-\ln G(\beta,x,z_{1},\ldots,z_{d})\rightarrow \exp\{-\sqrt{2d^{-1}\gamma}(x-\sum_{j=1}^{d}z_{j})\},
\end{eqnarray}
as $\beta\rightarrow0$. Similarly, we have
\begin{eqnarray}
\label{Le3.4.50}
-\ln G(\beta,y,z_{1},\ldots,z_{d})\rightarrow \exp\{-\sqrt{2d^{-1}\gamma}(y-\sum_{j=1}^{d}z_{j})\},
\end{eqnarray}
as $\beta\rightarrow0$.
Now, combining (\ref{Le3.4.30}-\ref{Le3.4.50}) with the dominated convergence theorem completes the proof of the case $r=d$.

\textbf{Lemma 3.4}. {\sl Under the conditions of Theorem 2.2, we have for $r=1, d$
\begin{eqnarray*}
\label{Le3.3.1}
&&\left|P\left(\Upsilon(O_{X}^{(r)})\geq -u_{\beta}(x), \tau(O_{X}^{(r)})\leq u_{\beta}(y)\right)\right.\\
&&\ \ \ \left. -P\left(O_{id}^{(r)}(X)\leq -i\beta, -u_{\beta}(x)-\beta^{-1}<i\leq -u_{\beta}(x), O_{kd}^{(r)}(X)\leq k\beta, u_{\beta}(y)\leq k< u_{\beta}(y)+\beta^{-1}  \right)\right|\rightarrow0
\end{eqnarray*}
as $\beta\rightarrow0$.
}

\textbf{Proof}.
 Obviously, the absolute value in Lemma 3.4 is bound above by
\begin{eqnarray*}
&&P\left(\bigcup_{i=-\infty}^{-u_{\beta}(x)-[\beta^{-1}]}(O_{id}^{(r)}(X)> -i\beta)\bigcup \bigcup_{k=u_{\beta}(y)+[\beta^{-1}]}^{\infty}(O_{kd}^{(r)}(X)> k\beta)\right)\\
&&\leq \sum_{i=-\infty}^{-u_{\beta}(x)-[\beta^{-1}]}P\left(O_{id}^{(r)}(X)> -i\beta\right)+\sum_{k=u_{\beta}(y)+[\beta^{-1}]}^{\infty}P(O_{kd}^{(r)}(X)> k\beta)\\
&&=\sum_{i=u_{\beta}(x)+[\beta^{-1}]}^{\infty}P\left(O_{id}^{(r)}(X)> i\beta\right)+\sum_{k=u_{\beta}(y)+[\beta^{-1}]}^{\infty}P(O_{kd}^{(r)}(X)> k\beta)\\
&&=\sum_{i=u_{\beta}(x)+[\beta^{-1}]}^{\infty}C_{d}^{r}(1-\Phi(i\beta))^{r}(1+o(1))
+\sum_{k=u_{\beta}(y)+[\beta^{-1}]}^{\infty}C_{d}^{r}(1-\Phi(k\beta))^{r}(1+o(1)),
\end{eqnarray*}
where in the last step, we use the fact that (\ref{Le3.2.2}). Using the fact (\ref{Le3.2.3}), the last sum is bounded above by
\begin{eqnarray*}
&&\sum_{i=u_{\beta}(x)+[\beta^{-1}]}^{\infty}\frac{C_{d}^{r}}{(2\pi)^{r/2}(i\beta)^{r}}\exp\left(-\frac{r(i\beta)^{2}}{2}\right)
+\sum_{k=u_{\beta}(y)+[\beta^{-1}]}^{\infty}\frac{C_{d}^{r}}{(2\pi)^{r/2}(k\beta)^{r}}\exp\left(-\frac{r(k\beta)^{2}}{2}\right)\\
&&\ll \beta^{-1} (\beta u_{\beta}(x))^{-(r+1)}\exp\left(-\frac{r(\beta u_{\beta}(x)+1)^{2}}{2}\right)+\beta^{-1}(\beta u_{\beta}(y))^{-(r+1)}\exp\left(-\frac{r(\beta u_{\beta}(y)+1)^{2}}{2}\right)\\
&&\ll\beta^{-\rho_{\beta}}|\ln \beta|^{-\frac{(1+r)\rho_{\beta}}{2}}e^{-K|\ln \beta|^{1/2}},
\end{eqnarray*}
which tends to 0, since $\rho_{\beta}\ln \beta^{-1}=\gamma$ under the conditions of Theorem 2.2.

\textbf{Proof of Theorem 2.2}. Define $Y_{ij}=\sqrt{1-\rho_{\beta}}\xi_{ij}+\sqrt{\rho_{\beta}}U_{j}$ and $O_{id}^{(r)}(Y)$  as in
Lemma 3.3.
Then by Lemmas 3.3 and 3.4, to prove Theorem 2.2, it suffices to prove for $r=1$ and $r=d$ that
\begin{eqnarray*}
&&\left|P\left(O_{id}^{(r)}(X)\leq -i\beta, -u_{\beta}(x)-\beta^{-1}<i\leq -u_{\beta}(x), O_{kd}^{(r)}(X)\leq k\beta, u_{\beta}(y)\leq k< u_{\beta}(y)+\beta^{-1}  \right)\right.\\
&&\left.-P\left(O_{id}^{(r)}(Y)\leq -i\beta, -u_{\beta}(x)-\beta^{-1}<i\leq -u_{\beta}(x), O_{kd}^{(r)}(Y)\leq k\beta, u_{\beta}(y)\leq k< u_{\beta}(y)+\beta^{-1}  \right)\right|\rightarrow0
\end{eqnarray*}
as $\beta\rightarrow0$. Using Lemma 3.1 again, we can bound the difference as
\begin{eqnarray*}
&&K\sum_{i=1}^{[\beta^{-1}]+1}|r_{i}-\rho_{\beta}|\sum_{j=u_{\beta}(x)}^{\infty}(j\beta)^{-2(r-1)}
\exp\left(-\frac{r((j\beta)^{2}+((j+i)\beta)^{2})}{2(1+\varpi_{i,\beta})}\right)\\
&&+K\sum_{i=1}^{[\beta^{-1}]+1}|r_{i}-\rho_{\beta}|\sum_{j=u_{\beta}(y)}^{\infty}(j\beta)^{-2(r-1)}
\exp\left(-\frac{r((j\beta)^{2}+((j+i)\beta)^{2})}{2(1+\varpi_{i,\beta})}\right)\\
&&+K \sum_{i=-u_{\beta}(x)-[\beta^{-1}]}^{-u_{\beta}(x)}\sum_{j=u_{\beta}(y)}^{u_{\beta}(y)+[\beta^{-1}]}|r_{|i-j|}-\rho_{\beta}|(-i\beta)^{-(r-1)}
(j\beta)^{-(r-1)}\exp\left(-\frac{r((i\beta)^{2}+(j\beta)^{2})}{2(1+\varpi_{|i-j|,\beta})}\right)\\
&&=:KR_{\beta}+KS_{\beta}+KT_{\beta},
\end{eqnarray*}
where $\varpi_{i,\beta}=\max\{|r_{i}|,\rho_{\beta}\}$.
Recall that $\delta=\sup_{i\geq1}|r_{i}|<1$ and $\delta_{n}=\sup_{i\geq n}|r_{i}|<1$. Choose $\varepsilon$ such that
$\varepsilon<(1-\delta)/(1+\delta)$ and $\varepsilon<1-|\ln \beta|^{-1/2}=:\theta$.  Split the above sum $J_{\beta}$ into three parts as
\begin{eqnarray*}
R_{\beta}&=&\sum_{i=1}^{[\beta^{-\varepsilon}]}|r_{i}-\rho_{\beta}|\sum_{j=u_{\beta}(x)}^{\infty}(j\beta)^{-2(r-1)}
\exp\left(-\frac{r((j\beta)^{2}+((j+i)\beta)^{2})}{2(1+\varpi_{i,\beta})}\right)\\
&+&\sum_{i=[\beta^{-\varepsilon}]+1}^{[\beta^{-\theta}]}|r_{i}-\rho_{\beta}|\sum_{j=u_{\beta}(x)}^{\infty}(j\beta)^{-2(r-1)}
\exp\left(-\frac{r((j\beta)^{2}+((j+i)\beta)^{2})}{2(1+\varpi_{i,\beta})}\right)\\
&+&\sum_{i=[\beta^{-\theta}]+1}^{[\beta^{-1}]+1}|r_{i}-\rho_{\beta}|\sum_{j=u_{\beta}(x)}^{\infty}(j\beta)^{-2(r-1)}
\exp\left(-\frac{r((j\beta)^{2}+((j+i)\beta)^{2})}{2(1+\varpi_{i,\beta})}\right)\\
&=:&R_{\beta,1}+R_{\beta,2}+R_{\beta,3}.
\end{eqnarray*}
By the same arguments as for $I_{\beta,1}$, we can show that $R_{\beta,1}=o(1)$ as $\beta\rightarrow0$. For the term $R_{\beta,2}$, let
$\pi_{\beta}=\sup_{i\geq [\beta^{-\theta}]}\varpi_{i,\beta}$. Then it is easy to see that $\pi_{\beta}|\ln\beta|\rightarrow\gamma$ as $\beta\rightarrow0$.
We have
\begin{eqnarray*}
R_{\beta,2}&\leq &(\beta^{-\theta}-\beta^{-\varepsilon})
\sum_{j=u_{\beta}(x)}^{\infty}(j\beta)^{-2(r-1)}\exp\left(-\frac{r(j\beta)^{2}}{1+\pi_{\beta}}\right)\\
&\ll&(\beta^{-\theta}-\beta^{-\varepsilon})\beta^{-1}(\beta u_{\beta}(x))^{-2r+1} \exp\left(-\frac{r(u_{\beta}(x)\beta)^{2}}{1+\pi_{\beta}}\right)\\
&\ll & \beta^{\frac{2(1-\rho_{\beta})}{1+\pi_{\beta}}-(1+\theta)}|\ln\beta|^{\frac{(1-\rho_{\beta})(1+r)}{1+\pi_{\beta}}-\frac{2r-1}{2}}\\
&\ll & \beta^{1-\theta}|\ln\beta|^{3/2}\\
&\ll& \exp\{-|\ln\beta|^{1/2}+3/2\ln|\ln \beta|\},
\end{eqnarray*}
which tends to 0 as as $\beta\rightarrow0$. To estimate the term $R_{\beta,3}$, using the following inequality
\begin{eqnarray}
\label{A0}
\left|r_{i}-\rho_{\beta}\right| \leqslant\left|r_{i}-\frac{\gamma}{\ln i}\right|+\gamma\left|\frac{1}{\ln i}-\frac{1}{\left|\ln \beta\right|}\right|
\end{eqnarray}
we have
\begin{eqnarray*}
R_{\beta,3}&\leq&\sum_{i=[\beta^{-\theta}]+1}^{[\beta^{-1}]+1}\left|r_{i}-\frac{\gamma}{\ln i}\right|\sum_{j=u_{\beta}(x)}^{\infty}(j\beta)^{-2(r-1)}
\exp\left(-\frac{r((j\beta)^{2}+((j+i)\beta)^{2})}{2(1+\varpi_{i,\beta})}\right)\\
&+&\gamma\sum_{i=[\beta^{-\theta}]+1}^{[\beta^{-1}]+1}\left|\frac{1}{\ln i}-\frac{1}{\left|\ln \beta\right|}\right|\sum_{j=u_{\beta}(x)}^{\infty}(j\beta)^{-2(r-1)}
\exp\left(-\frac{r((j\beta)^{2}+((j+i)\beta)^{2})}{2(1+\varpi_{i,\beta})}\right)\\
&=:&R_{\beta,31}+R_{\beta,32}.
\end{eqnarray*}
 By the same arguments as for the proof of $I_{\beta,2}$ in Theorem 2.1, we have
\begin{eqnarray*}
R_{\beta,31}&\ll& |\ln \beta|^{-1}\sum_{i=[\beta^{-\theta}]+1}^{[\beta^{-1}]+1}\left|r_{i}\ln i-\gamma\right|\sum_{j=u_{\beta}(x)}^{\infty}(j\beta)^{-2(d-1)}
\exp\left(-\frac{r((j\beta)^{2}+((j+i)\beta)^{2})}{2(1+\varpi_{i,\beta})}\right)\\
&\leq&\sup_{i\geq \beta^{-\theta}}\{|r_{i}\ln i-\gamma|\}\beta^{-2(\pi_{\beta}+\rho_{\beta})}|\ln\beta|^{-(1+r)(\rho_{\beta}+\pi_{\beta})},
\end{eqnarray*}
which tends to 0, since $r_{i}\ln i\rightarrow\gamma$.
Using the inequality for $[\beta^{-\theta}]+1 \leq i\leq[\beta^{-1}]+1$
$$
\left|\frac{1}{\ln i}-\frac{1}{\left|\ln \beta\right|}\right|\ll|\ln \beta|^{-3/2}
$$
and by the same arguments as for the proof of $R_{\beta,31}$, we have
$$R_{\beta,32}\leq \beta^{-2(\pi_{\beta}+\rho_{\beta})}|\ln\beta|^{-(1+r)(\pi_{\beta}+\rho_{\beta})-1/2},$$
which also tends to 0.
By the same arguments, we get $S_{\beta}\rightarrow0$  as $\beta\rightarrow0$.\\
For the term $T_{\beta}$, Note that $u_{\beta}(x)\sim u_{\beta}(y)\sim u_{\beta}:=\beta^{-1}(\frac{2}{r}\ln \beta^{-1})^{1/2}$ as $\beta\rightarrow0$.
Define
$\kappa_{\beta}=\sup_{i\geq 2[u_{\beta}]}\varpi_{i,\beta}$.
Using (\ref{A0}) again, we have
\begin{eqnarray*}
K_{\beta}&\leq&
\sum_{i=-u_{\beta}(x)-[\beta^{-1}]}^{-u_{\beta}(x)} \sum_{j=u_{\beta}(y)}^{u_{\beta}(y)+[\beta^{-1}]}\left|r_{|i-j|}-\frac{\gamma}{\ln|i-j|}\right|(-i\beta)^{-(r-1)}(j\beta)^{-(r-1)}
\exp\left(-\frac{r(i\beta)^{2}+r(j\beta)^{2}}{2(1+\kappa_{\beta})}\right)\\
&+&\gamma \sum_{i=-u_{\beta}(x)-[\beta^{-1}]}^{-u_{\beta}(x)} \sum_{j=u_{\beta}(y)}^{u_{\beta}(y)+[\beta^{-1}]}\left|\frac{1}{\ln|i-j|}-\frac{1}{|\ln \beta|}\right|(-i\beta)^{-(r-1)}(j\beta)^{-(r-1)}
\exp\left(-\frac{r(i\beta)^{2}+r(j\beta)^{2}}{2(1+\kappa_{\beta})}\right)\\
&=:& K_{\beta,1}+K_{\beta,2}
\end{eqnarray*}
 By the same arguments as for the proof of $R_{\beta,31}$, we have
\begin{eqnarray*}
 K_{\beta,1}&\leq&\sup_{i\geq 2[u_{\beta}]}\{|r_{i}\ln i-\gamma|\}|\ln \beta|^{-1} \sum_{i=u_{\beta}(x)}^{u_{\beta}(x)+[\beta^{-1}]}(i\beta)^{-(r-1)}
\exp\left(-\frac{r(i\beta)^{2}}{2(1+\kappa_{\beta})}\right)\\
\ \ \ \ \ \ \ \ &\times&\sum_{j=u_{\beta}(y)}^{u_{\beta}(y)+[\beta^{-1}]}(j\beta)^{-(r-1)}
\exp\left(-\frac{r(j\beta)^{2}}{2(1+\kappa_{\beta})}\right)\\
&\ll&\sup_{i\geq 2[u_{\beta}]}\{|r_{i}\ln i-\gamma|\}\beta^{-2(\kappa_{\beta}+\rho_{\beta})}|\ln\beta|^{-(1+r)(\rho_{\beta}+\kappa_{\beta})},
\end{eqnarray*}
which tends to 0, since $r_{n}\ln n\rightarrow\gamma$ as $n\rightarrow\infty$.
Using the inequality for $u_{\beta}(x)+u_{\beta}(y)<|i-j|<u_{\beta}(x)+u_{\beta}(y)+2[\beta^{-1}]$
$$
\left|\frac{1}{\ln |i-j|}-\frac{1}{\left|\ln \beta\right|}\right|\ll|\ln \beta|^{-3/2},
$$
we have
\begin{eqnarray*}
 K_{\beta,2}&\leq&|\ln \beta|^{-3/2}\sum_{i=u_{\beta}(x)}^{u_{\beta}(x)+[\beta^{-1}]}(i\beta)^{-(r-1)}
\exp\left(-\frac{r(i\beta)^{2}}{2(1+\kappa_{\beta})}\right)\\
&\times&\sum_{j=u_{\beta}(y)}^{u_{\beta}(y)+[\beta^{-1}]}(j\beta)^{-(r-1)}
\exp\left(-\frac{r(j\beta)^{2}}{2(1+\kappa_{\beta})}\right)\\
&\ll&\beta^{-2(\kappa_{\beta}+\rho_{\beta})}|\ln\beta|^{-(1+r)(\rho_{\beta}+\kappa_{\beta})-1/2},
\end{eqnarray*}
which tends to 0 as $\beta\rightarrow0$.
 Thus, we have $K_{\beta}\rightarrow0$  as $\beta\rightarrow0$.
The proof of Theorem 2.2 is complete.

\subsection{Proof of Theorem 2.3}

We need the following two lemmas to prove Theorem 2.3.

\textbf{Lemma 3.5}. {\sl Let $\rho_{\beta}$ be  a function of $\beta$ with $0\leq \rho_{\beta}<c<1$ and
$Y_{\beta,i}=\sqrt{1-\rho_{\beta}}\xi_{i}+\sqrt{\rho_{\beta}}U$, where $\xi_{i}$ is a sequence of independent standard Gaussian
variables  and $U$ is a standard Gaussian random variable independent with $\xi_{i}$.
 If $\rho_{\beta}|\ln \beta|\rightarrow\infty$ as $\beta\rightarrow0$, we have for any $x,y\in \mathbb{R}$
\begin{eqnarray*}
\label{Le3.5.1}
\lim_{\beta\rightarrow0}P\left(Y_{\beta,i}\leq -i\beta, -u_{\beta}(x)-\beta^{-1}<i\leq -u_{\beta}(x),
Y_{\beta,j}\leq j\beta, u_{\beta}(y)\leq j< u_{\beta}(y)+\beta^{-1} \right)=\Phi(\min(x,y)),
\end{eqnarray*}
where $u_{\beta}(x)$ is defined in (\ref{Th2.2.2}) with $\rho_{\beta}$ satisfying the above assumptions.
}

\textbf{Proof}. Since
\begin{eqnarray*}
&&P\left(Y_{\beta,i}\leq -i\beta, -u_{\beta}(x)-\beta^{-1}<i\leq -u_{\beta}(x),
Y_{\beta,j}\leq j\beta, u_{\beta}(y)\leq j< u_{\beta}(y)+\beta^{-1} \right)\\
&&=P\left(\sqrt{1-\rho_{\beta}}\xi_{i}+\sqrt{\rho_{\beta}}U\leq -i\beta, -u_{\beta}(x)-\beta^{-1}<i\leq -u_{\beta}(x),\right.\\
&&\ \ \ \ \ \ \ \ \ \ \ \left.\sqrt{1-\rho_{\beta}}\xi_{j}+\sqrt{\rho_{\beta}}U\leq j\beta, u_{\beta}(y)\leq j< u_{\beta}(y)+\beta^{-1} \right)\\
&&=\int_{-\infty}^{+\infty}
P\left(\xi_{i}\leq -\frac{\beta(i+\sqrt{\rho_{\beta}}z/\beta)}{\sqrt{1-\rho_{\beta}}},  -u_{\beta}(x)-\beta^{-1}<i\leq -u_{\beta}(x),\right.\\
&&\ \ \ \ \ \ \ \ \ \ \ \ \ \ \ \ \left.  \xi_{j}\leq -\frac{\beta(j+\sqrt{\rho_{\beta}}z/\beta)}{\sqrt{1-\rho_{\beta}}},  u_{\beta}(y)\leq j< u_{\beta}(y)+\beta^{-1} \right)d\Phi(z),
\end{eqnarray*}
by the definition of $\xi_{i}$, we get  the above integral equals
\begin{eqnarray}
\label{Le3.5.1}
&&\int_{-\infty}^{+\infty}
P\left(\xi_{i}\leq -\frac{i\beta}{\sqrt{1-\rho_{\beta}}}, -(u_{\beta}(x)-\sqrt{\rho_{\beta}}z/\beta)-\beta^{-1}<i\leq -(u_{\beta}(x)-\sqrt{\rho_{\beta}}z/\beta)\right)\nonumber\\
&&\ \ \ \ \ \ \times P\left(\xi_{j}\leq \frac{j\beta}{\sqrt{1-\rho_{\beta}}},  (u_{\beta}(y)-\sqrt{\rho_{\beta}}z/\beta)\leq j <(u_{\beta}(y)-\sqrt{\rho_{\beta}}z/\beta)+\beta^{-1}\right)d\Phi(z)\nonumber\\
&&=:\int_{-\infty}^{+\infty}F(\beta,x,z)F(\beta,y,z)d\Phi(z).
\end{eqnarray}
By the same arguments as in Lemma 3.3 (see also the proof of Lemma 2.2 of H\"{u}sler (1977)), we have
$$-\ln F(\beta,x,z)\rightarrow \Bigg\{
\begin{array}{l}
0, \ \mbox{if}\ \ x>z;\\
\infty, \ \mbox{if}\ \ x<z,
\end{array} $$
and
$$-\ln F(\beta,y,z)\rightarrow \Bigg\{
\begin{array}{l}
0, \ \mbox{if}\ \ y>z;\\
\infty, \ \mbox{if}\ \ y<z,
\end{array} $$
which combining with (\ref{Le3.5.1}) completes the proof of the lemma.

\textbf{Lemma 3.6}. {\sl Under the conditions of Theorem 2.3, we have
\begin{eqnarray*}
\label{Le3.6.1}
&&\left|P\left(\Upsilon(X)\geq -u_{\beta}(x), \tau(X)\leq u_{\beta}(x)\right)\right.\\
&&\ \ \ \ \ \ \ \ \left.-P\left(X_{i}\leq -i\beta, -u_{\beta}(x)-[\beta^{-1}]<i\leq -u_{\beta}(x), X_{j}\leq j\beta, u_{\beta}(y)\leq i< u_{\beta}(y)+[\beta^{-1}]
 \right)\right|\rightarrow0
\end{eqnarray*}
as $\beta\rightarrow0$.
}

\textbf{Proof}. Let $A_{X}=\left(\Upsilon(X)\geq -u_{\beta}(x)\right)$,
$A'_{X}=\left(X_{i}\leq -i\beta, -u_{\beta}(x)-[\beta^{-1}]<i\leq -u_{\beta}(x)\right)$, $B_{X}=\left(\tau(X)\leq u_{\beta}(x)\right)$ and
$B'_{X}=\left(X_{j}\leq j\beta, u_{\beta}(y)\leq i< u_{\beta}(y)+[\beta^{-1}]\right)$. Obviously,
\begin{eqnarray*}
\left|P(A_{X}\bigcap B_{X})-P(A'_{X}\bigcap B'_{X})\right|&\leq&\left|P(A_{X}\bigcap B_{X})-P(A'_{X}\bigcap B_{X})\right|\\
&+&\left|P(A'_{X}\bigcap B_{X})-P(A'_{X}\bigcap B'_{X})\right|\\
&\leq &\left|P(A_{X})-P(A'_{X})|+|P(B_{X})-P(B'_{X})\right|\\
&\leq & P(E_{X}^{c})+P(F_{X}^{c}),
\end{eqnarray*}
where
$E_{X}=\left(\Upsilon(X)\geq -u_{\beta}(x)-[\beta^{-1}]\right)$
and
$F_{X}=\left(\tau(X)\leq u_{\beta}(x)+[\beta^{-1}]\right)$.
From the proof of Lemma 2.3 of H\"{u}sler (1977),
we get $P(E_{X})\rightarrow1$, as $\beta\rightarrow0$.
Similarly, we have
$P(F_{X})\rightarrow1$, as $\beta\rightarrow0$.
The proof of the lemma is complete.

\textbf{Proof of Theorem 2.3}. By Lemma 3.6, it suffices to show that
\begin{eqnarray*}
P\left(X_{i}\leq -i\beta, -u_{\beta}(x)-[\beta^{-1}]<i\leq -u_{\beta}(x), X_{j}\leq j\beta, u_{\beta}(y)\leq i< u_{\beta}(y)+[\beta^{-1}]
 \right)\rightarrow\Phi(\min(x,y))
\end{eqnarray*}
as $\beta\rightarrow0$.
Define $Y_{\beta,i}=\sqrt{1-\rho_{\beta}}\xi_{i}+\sqrt{\rho_{\beta}}U$ as in
Lemma 3.5. Since $r_{i}\geq \rho_{\beta}$ for all $i\leq \beta^{-1}$, by Slepian's lemma (see, e.g. Leadbetter et al. (1983)),
we have
\begin{eqnarray*}
&&\lim_{\beta\rightarrow0}P\left(X_{i}\leq -i\beta, -u_{\beta}(x)-[\beta^{-1}]<i\leq -u_{\beta}(x), X_{j}\leq j\beta, u_{\beta}(y)\leq i< u_{\beta}(y)+[\beta^{-1}]
 \right)\\
&&\geq \lim_{\beta\rightarrow0}P\left(Y_{\beta,i}\leq -i\beta, -u_{\beta}(x)-\beta^{-1}<i\leq -u_{\beta}(x),
Y_{\beta,j}\leq j\beta, u_{\beta}(y)\leq j< u_{\beta}(y)+\beta^{-1} \right),
\end{eqnarray*}
which by Lemma 3.5 equals $\Phi(\min(x,y))$. Thus, we only need to show
\begin{eqnarray*}
&&P\left(X_{i}\leq -i\beta, -u_{\beta}(x)-[\beta^{-1}]<i\leq -u_{\beta}(x), X_{j}\leq j\beta, u_{\beta}(y)\leq i< u_{\beta}(y)+[\beta^{-1}]
 \right)\\
 &&\ \ \ \ \leq \Phi(\min(x,y)+\varepsilon)
 \end{eqnarray*}
for all $\varepsilon>0$. Since $r_{n}$ is convex, there is a Gaussian sequence $\{Z_{i}=Z_{i}(\beta), i\in \mathbb{Z}\}$ with the correlations
(see, e.g. Mittal and Ylvisaker (1975))
$$\varrho_{k}=(r_{k}-\rho_{\beta})/(1-\rho_{\beta})\ \ \mbox{for}\ \ k=1,2,\ldots, [\beta^{-1}].$$
Let $V$ be independent of $\{Z_{i}, i\in \mathbb{Z}\}$, such that
$$X_{i}=\sqrt{1-\rho_{\beta}}Z_{i}+\sqrt{\rho_{\beta}}V.$$
Now
\begin{eqnarray*}
&&P\left(X_{i}\leq -i\beta, -u_{\beta}(x)-[\beta^{-1}]<i\leq -u_{\beta}(x), X_{j}\leq j\beta, u_{\beta}(y)\leq i< u_{\beta}(y)+[\beta^{-1}]
 \right)\\
&&=\int_{-\infty}^{\infty}P\left(Z_{i}\leq -\frac{i\beta}{\sqrt{1-\rho_{\beta}}}, -u_{\beta}(x-z)-[\beta^{-1}]<i\leq -u_{\beta}(x-z),\right.\\
&&\ \ \ \ \ \ \ \ \ \ \ \ \  \left. Z_{j}\leq \frac{j\beta}{\sqrt{1-\rho_{\beta}}}, u_{\beta}(y-z)\leq j< u_{\beta}(y-z)+[\beta^{-1}] \right)d\Phi(z) \\
&&\leq \Phi\left(\min(x,y)+\varepsilon\right)+\int_{\min(x,y)+\varepsilon}^{\infty}P\left(Z_{i}\leq -\frac{i\beta}{\sqrt{1-\rho_{\beta}}}, -u_{\beta}(x-z)-[\beta^{-1}]<i\leq -u_{\beta}(x-z),\right.\\
&&\ \ \ \ \ \ \ \ \ \ \ \ \  \left. Z_{j}\leq \frac{j\beta}{\sqrt{1-\rho_{\beta}}}, u_{\beta}(y-z)\leq j< u_{\beta}(y-z)+[\beta^{-1}] \right)d\Phi(z)\\
&&\leq\Phi\left(\min(x,y)+\varepsilon\right)+P\left(Z_{i}\leq -\frac{i\beta}{\sqrt{1-\rho_{\beta}}}, -u_{\beta}(-\varepsilon)-[\beta^{-1}]<i\leq -u_{\beta}(-\varepsilon)\right).
\end{eqnarray*}
Therefore, to complete the proof of the theorem, it suffices to show that the second term tends to $0$,
which has been done in the proof of Theorem 2.3. in H\"{u}sler (1977).

\bigskip

{\bf Acknowledgement}: The authors would like to thank the two
referees for the valuable suggestions.

\end{document}